\documentclass[notitlepage,11pt]{article}
\usepackage{amssymb,amsfonts,amsmath,geometry,esint
}

\catcode`\@=11
\@addtoreset{equation}{section}

\catcode`\@=12

\usepackage{color}

\def\R{\mathbb{R}}

\def\tu{\widetilde u}

\def\f{\varphi}

\def\irn{\int\limits_{\R^n}}
\def\iirn{\int\limits_0^\infty\!\!\int\limits_{\R^n}}

\def\fintn{\fint\limits_{\S^{n-1}}}
\def\eps{\varepsilon}
\def\div{{\rm div}}

\def\Ds{\left(-\Delta\right)^{\!s}} 
\def\Dshalf{\left(-\Delta\right)^{\!\frac{s}{2}}}  

\def\S{\mathbb S}

\def\sstar{{2^*_s}}

\def\proof{\noindent{\textbf{Proof. }}}
\def\QED{\hfill {$\square$}\goodbreak \medskip}

\newtheorem{Theorem}{Theorem}[section]
\newtheorem{Lemma}[Theorem]{Lemma}

\newtheorem{Corollary}[Theorem]{Corollary}
\newtheorem{Remark}[Theorem]{Remark}

\usepackage{mathtools}

\begin{document}

\title 
{A tool for symmetry breaking and multiplicity \\in some nonlocal problems}

\author{Roberta Musina
\footnote{Dipartimento di Matematica ed Informatica, Universit\`a di Udine,
via delle Scienze, 206 -- 33100 Udine, Italy. Email: {roberta.musina@uniud.it}. 
}~ 
and \setcounter{footnote}{6}
Alexander I. Nazarov
\footnote{
St.Petersburg Dept of Steklov Institute, Fontanka 27, St.Petersburg, 191023, Russia, 
and St.Petersburg State University, 
Universitetskii pr. 28, St.Petersburg, 198504, Russia. E-mail: al.il.nazarov@gmail.com.
}
}

\date{}

\maketitle


\noindent
{\small {\bf Abstract.} We prove some basic inequalities relating the Gagliardo-Nirenberg seminorms of a symmetric function
$u$ on $\R^n$ and of its perturbation $u\f_\mu$, where $\f_\mu$ is a suitably chosen eigenfunction of the
Laplace-Beltrami operator on the sphere $\S^{n-1}$, thus providing a technical but rather powerful tool to detect
symmetry breaking and multiplicity phenomena in variational equations driven by the fractional Laplace operator.
A concrete application to a problem related to the fractional Caffarelli-Kohn-Nirenberg inequality is given.}


\medskip\noindent
{\small {\bf 2010 Mathematics Subject Classfication:}}  35R11; 35B06; 35A02.

\normalsize

\bigskip

\section{Introduction} 

Let
${\f_{\!\mu}}$ be an eigenfunction of the Laplace-Beltrami
operator on the sphere $\S^{n-1}$, $n\ge 2$, relative to the positive eigenvalue $\mu$ and normalized by the condition
\begin{equation}
\label{eq:normalization}
\fintn {\f^2_{\!\mu}}~\!d\sigma=1.
\end{equation}
Take a  ''good'' radial function $u$ on $\R^n$. We break the  symmetry of $u$
by defining, via polar coordinates, $(u{\f_{\!\mu}})(r\sigma)=u(r){\f_{\!\mu}}(\sigma)$.
The function $u{\f_{\!\mu}}$ has the same $L^2$ norm as $u$ and it is orthogonal to $u$
in $L^2$.
A popular and efficient technique to detect 
symmetry breaking and multiplicity results for a large class of variational problems
is based on the comparison between the energies of $u$ and $u{\f_{\!\mu}}$. We cite for instance \cite{SSW, GM2}, where the trivial equality
\begin{equation}
\label{eq:trivial}
\irn|\nabla (u{\f_{\!\mu}})|^2~\!dx-\irn|\nabla u|^2~\!dx
= \mu\irn|x|^{-2}|u|^2~\!dx
\end{equation}
is crucially used to tackle certain problems driven by the Laplace operator $-\Delta$. 
We cite also  \cite{N04, Sch, KN, MLincei, CMMilan, NN}, 
where the $p$-Laplacian or more general second order, possibly degenerate operators in divergence form 
are considered, and
\cite{BM, CC}, that deal with fourth-order variational equations.

In dealing with variational problems involving the Dirichlet Laplacian $\Ds$, $0<s<1$,
a simple and powerful identity comparable with (\ref{eq:trivial}) is hopeless.
However, as a corollary of the more general Lemma \ref{L:Lbs} in Section \ref{S:perturbing}, we obtain the existence
of a positive constant $c_{\mu}$, not depending on $u$, such that 
\begin{equation}
\label{eq:uw_int}
\irn|\Dshalf (u{\f_{\!\mu}})|^2~\!dx-\irn|\Dshalf u|^2~\!dx \le c_{\mu}\irn |x|^{-2s}|u|^2~\!dx~\!.
\end{equation}
Besides its impact on the study of fractional differential equations, 
inequality (\ref{eq:uw_int}) might have an independent interest. For instance, it is strongly related to
Bochner's relations and to the results in \cite{DKK, GT}. 

Inequality (\ref{eq:uw_int}) and its generalizations below provide quite useful technical tools. 
In order to illustrate their applications in concrete
problems
we take as model the fractional Hardy-Sobolev  inequality 
\begin{equation}
\label{eq:HS_local}
\irn|\Dshalf u|^2~\!dx+\lambda\irn|x|^{-2s}|u|^2~\!dx\ge \mathcal S^{\lambda}_q\Big(\irn |x|^{-bq}|u|^q~\!dx\Big)^{\frac2q}~\!,
\quad u\in \mathcal D^s(\R^n)
\end{equation}
and its associated semilinear  Euler-Lagrange equation 
\begin{equation}
\label{eq:CKN}
\Ds u+\lambda|x|^{-2s} u=|x|^{-bq}u^{q-1}~,\quad  u\in\mathcal D^s(\R^n)~\!.
\end{equation}
Inequality (\ref{eq:HS_local}) is clearly related
to the fractional Hardy and Sobolev inequalities,
$$
\irn|\Dshalf u|^2~\!dx\ge H_{\!s}\irn |x|^{-2s}|u|^2~\!dx~\!,\quad
\irn|\Dshalf u|^2~\!dx\ge S_{\sstar}\Big(\irn |u|^\sstar~\!dx\Big)^{\frac2\sstar}~\!,
$$
where 
$\sstar=\frac{2n}{n-2s}$
is the critical Sobolev exponent.
The explicit values of the  Hardy constant $H_{\!s}$ and of the Sobolev constant  $S_\sstar$  have been computed in \cite{He}, \cite{CoTa}, respectively.

From now on, we take exponents $q, b$ satisfying
\begin{equation}
2<q<\sstar=\frac{2n}{n-2s}~,\qquad 
\frac{n}{q}-b=\frac{n}{2}-s~\!.
\label{eq:bq}
\end{equation}
By using H\"older's interpolation, it is easy to see that 
(\ref{eq:HS_local}) holds with a positive best constant $S^\lambda_q$ provided that 
$\lambda>-H_{\!s}$. Under these assumptions, nowadays standard 
arguments give the existence of an extremal for the noncompact minimization problem
\begin{equation}
\label{eq:min_problem}
S^\lambda_q=\inf_{u\in\mathcal D^s(\R^n)\atop u\neq 0} J_{\lambda}(u)~,\qquad
J_{\lambda}(u):=\frac{\|\Dshalf u\|_2^2+\lambda\||x|^{-s} u\|_2^2}{\||x|^{-b} u\|_q^2}~\!,
\end{equation}
see Corollary \ref{C:existence} in Section \ref{S:variational}. 
Thanks to
(\ref{eq:uw_int}), in Section \ref{S:proofs} we  prove the next result.

\begin{Theorem}
\label{T:BS}
If $\lambda>0$ is large enough then symmetry breaking occurs, that is,
no extremal for $S^{\lambda}_q$ is radially symmetric.
\end{Theorem}

Up to a Lagrange multiplier, any extremal $\underline{u}^\lambda$ 
for $S^{\lambda}_q$ is nonnegative and solves (\ref{eq:CKN}) in the weak sense. By Theorem \ref{T:BS} and Lemma \ref{L:invariance} in
Section \ref{S:variational},
we have that  for $\lambda$ large enough, problem (\ref{eq:CKN}) has in fact two distinct 
nonnegative solutions: the ground state solution $\underline u^{\lambda}$, which is not radial, 
and a radial one, 
that minimizes $J_{\lambda}(u)$ on the space of radial functions. 

The exploitation of different symmetries and Lemma \ref{L:Lbs} 
lead to the next multiplicity results.

\begin{Theorem}
\label{T:hradial}
Let $n\ge 2$ be even. For any integer $h\ge1$, there exists $\lambda_h>0$ such that
for $\lambda>\lambda_h$, problem (\ref{eq:CKN}) has at least $h$  nonnegative solutions, that 
are distinct modulo rotations of $\R^n$.
\end{Theorem}

\begin{Theorem}
\label{T:dimension}
Let $n\ge 3$ be odd. Then for $\lambda$ large enough problem (\ref{eq:CKN}) has $N(n)\ge4$
nonnegative solutions, that are distinct modulo rotations of $\R^n$.
\end{Theorem}

Let us conclude this introduction by pointing out few facts.

If $\lambda\in(-H_s,0]$, then any nonnegative solution to (\ref{eq:CKN}) is radially symmetric about the origin;
for the proof, notice that $b>0$ and adapt the moving plane argument in \cite{DMPS}. In particular, any extremal for $S^{\lambda}_q$ is radially symmetric if $\lambda\le 0$ (the last statement 
can be also proved by the Schwarz symmetrization, see, e.g., \cite[Theorem 2.31 and Sec. II.2]{Kawohl}).

By Theorem \ref{T:BS}, there exists an optimal parameter $\widehat\lambda=\widehat\lambda(n,s,q)\ge  0$ such that symmetry breaking occurs whenever $\lambda>\widehat\lambda$. 
In the local case $s=1$, symmetry breaking has been firstly pointed out in \cite{CaWa}.\footnote{The equation considered in \cite{CaWa} is related to (\ref{eq:CKN}) with $s=1$ 
via functional change $u(x)\mapsto |x|^{a}u(x)$.}
Nowadays the picture is complete, thanks to the results in \cite{FeSc, DET, DEL}, one gets that $\widehat\lambda=\frac{4(n-1)}{q^2-4}$.

In the nonlocal case we cannot prove even whether the set of $\lambda$ providing the symmetry breaking is connected, nor we have any conjecture about the value of $\widehat\lambda$.
By exploiting our proof  and thanks to Remark \ref{R:c}, 
one can find rough upper bounds on $\widehat\lambda$
in case $n\ge 4$, or $n=3$ and $0<s\le \frac12$.

%

The paper is organized as follows. In the next two sections we prove some crucial inequalities, including (\ref{eq:uw_int}).
The main tools are the Caffarelli-Silvestre extension technique \cite{CS} and the results in \cite{MNine}. Section \ref{S:variational} contains the main variational tools
and a criterion to distinguish solutions to (\ref{eq:CKN}) enjoing different symmetry properties, see the {\em Basic Lemma} \ref{L:Cnk}. The proofs of the main Theorems
are collected in Section \ref{S:proofs}. 

\begin{Remark}
Minor modifications in the variational arguments give simmetry breaking and multiplicity of positive solutions to the Dirichlet problem
$$
\begin{cases}
\Ds u= u^{q-1}&\text{in $A= \{R<|x|<R+1\}\subset\R^n$}\\
u=0& \text{in $\R^n\setminus A$}
\end{cases}
$$ 
for $q\in(2,\sstar)$ and $R$ large. See however \cite{Us_He}, where a different argument is used.
\end{Remark}

\section{Preliminaries}
\label{S:preliminaries}

The fractional Laplacian $\Ds$ in $\R^n$, $n\ge 2$, is formally defined by 
$$
{\mathcal F}\big[\Ds u\big] = |\xi|^{2s}{\mathcal F}[u]~\!,
$$
where ${\mathcal F}=\displaystyle{{\mathcal F}[u](\xi)= (2\pi)^{-\frac{n}{2}}
\int_{\R^n} e^{-i~\!\!\xi\cdot x}u(x)~\!dx}$ is the Fourier transform.

\medskip

Thanks to the Sobolev inequality, the space
$$
\mathcal D^s(\R^n)=\big\{u\in L^\sstar(\R^n)~|~\Dshalf u \in L^2(\R^n)~\!\big\}
$$
naturally inherits a Hilbertian structure from the scalar product
$$
(u,v)=\irn\Dshalf u\Dshalf v~\!dx=\irn|\xi|^{2s}\mathcal F[u]\overline{\mathcal F[v]}~\!d\xi~\!.
$$
From now on, we will always use the shorter notation $\mathcal D^s$ instead of $\mathcal D^s(\R^n)$. 

\medskip

In the breakthrough paper \cite{CS}, Caffarelli and Silvestre 
investigated the  relations between the nonlocal operator $\Ds$ in
$\mathbb R^n\ni x$ and the pointwise defined differential operator $-\text{div}(y^{1-2s}\nabla)$ in
$\R^{n+1}_+\equiv \R^n\times(0,\infty)\ni (x,y)$. It turns out that 
any function $w$ in the space
$$
\mathcal W^s=\mathcal W^s(\R^{n+1}_+)=\big\{w:\R^{n+1}_+\to\R~~\text{measurable,}~~
\int\limits_0^\infty\!\int\limits_{\mathbb{R}^n} y^{1-2s }|\nabla w|^2\,dxdy<\infty~\big\}
$$
has a trace on the boundary of $\R^{n+1}_+$,  $w|_{\partial \R^{n+1}_+}\in \mathcal D^s$
and for any $u\in\mathcal D^s$ we have that 
\begin{equation}
\label{eq:CS_problem}
\irn|\Dshalf u|^2~\!dx=\inf_{{w\in\mathcal W^s \atop w|_{\partial\R^{n+1}_+}=u}}
C_s
\iirn y^{1-2s }|\nabla w|^2\,dxdy~\!,\quad C_s=\frac {\Gamma(s)}{2^{1-2s}\Gamma(1-s)}.
\end{equation}
The {\em Caffarelli-Silvestre extension} $w_u$ of a function $u\in\mathcal D^s$ is the unique solution to the convex minimization problem in (\ref{eq:CS_problem}), hence it satisfies
\begin{gather}
\nonumber
-\div\big( y^{1-2s}\nabla w_u\big)=0\quad \text{in $\R^{n+1}_+$,}\quad w_u(x,0)=u(x),\\
C_s~\! \iirn y^{1-2s }|\nabla w_u|^2\,dxdy=\irn|\Dshalf u|^2~\!dx~\!.
\label{quad_D}
\end{gather}

Recall that the Hardy type inequality
$$
\int\limits_0^\infty\!\irn y^{1-2s}|\nabla w|^2\,dxdy\ge \Big(\frac {n-2s}{2}\Big)^2 \int\limits_0^\infty\!\int\limits_{\mathbb{R}^n} y^{1-2s }\,\frac {|w|^2}{|x|^2+y^2}\,dxdy
~,
\quad w\in \mathcal W^s
$$
holds with a sharp and not achieved constant,
see \cite[Section 2]{MNaihp}. In particular, $\mathcal W^s$  inherits a natural 
Hilbert space structure and the map $u\mapsto w_u$ is an  isometry,
up to the constant $C_s$. 
In the next lemma we provide a crucial relation between the Hardy integrals of $u\in \mathcal D^s$ and 
of its extension $w_u\in\mathcal W^s$.

\begin{Lemma} 
\label{L:crucial}
Let  $w_u$ be the Caffarelli-Silvestre extension of $u\in \mathcal D^s$. Then
\begin{equation}
\label{eq:uw0}
\iirn y^{1-2s}\frac{|w_u|^2}{|x|^2+y^2}~\!dxdy 
\le \gamma \irn\frac{|u|^2}{|x|^{2s}}~\!dx~\!,
\end{equation}
where the positive constant $\gamma$ does not depend on $u$.
\end{Lemma}

\proof
By \cite[Theorem 1]{MNine} there exists  a constant $\hat{c}>0$, possibly  depending  on $n, s$ but not on $u$, such that  
\begin{equation}
\label{eq:dis}
\irn \frac{|w_u(x,y)|^2}{|x|^2+y^2}~\!dx\le \hat{c}\irn \frac{|u(x)|^2}{|x|^2+y^2}~\!dx\quad\text{for any $y>0$.}
\end{equation}
Since 
$$
\iirn y^{1-2s}\frac{|u(x)|^2}{|x|^2+y^2}~\!dxdy=\irn |u(x)|^2~dx\int\limits_0^\infty\frac{y^{1-2s}}{|x|^2+y^2}~\!dy=
\frac12 \Gamma(s)\Gamma(1-s) \irn \frac{|u(x)|^2}{|x|^{2s}}~\!dx,
$$
the conclusion follows immediately, with $\gamma=\hat{c}~\!\Gamma(s)\Gamma(1-s)/2$.
\QED

\begin{Remark}
\label{R:c}
By \cite[Theorem 2]{MNine} we know that (\ref{eq:dis}) holds with $\hat{c}= 1$, provided that
$n\ge 4$ or $n=3$ and $0<s\le \frac12$. In this case we obtain the estimate
\begin{equation}
\label{eq:uw1}
\iirn y^{1-2s}\frac{|w_u|^2}{|x|^2+y^2}~\!dxdy 
\le \frac{\pi}{2\sin(\pi s)} \irn\frac{|u|^2}{|x|^{2s}} ~\!dx\quad\text{for any $u\in\mathcal D^s$}.
\end{equation}
We conjecture that (\ref{eq:uw1}) holds with a sharp constant, at least for $n\ge 3$. The lowest dimensional case $n=2$
looks more obscure. Finally, it would be of interest to investigate whether (\ref{eq:uw0}) holds in case $n=1$, $s\in(0,\frac12)$.
\end{Remark}

\section{Perturbing symmetric functions} 
\label{S:perturbing}

Let $n=km$ with $k\ge 2$, $m\ge 1$, and write $\R^n$ as the Cartesian product of $m$ copies of $\R^k$.
It is convenient to denote by $\R^k_j$ the $j$-th copy of $\R^k$, so that $\R^n=\R^k_1\times\dots\times\R^k_m$.
The variable in $\R^k_j$ is $x_j$;  its polar coordinates are $r_j=|x_j|$, $\sigma_j\in \S^{k-1}_j$, where $\S^{k-1}_j$ is the unit sphere
in $\R^k_j$.

\medskip

In the next crucial lemma we take 
a proper closed subgroup 
$\mathcal G_k$  of $O(k)$ and an eigenfunction $\phi$ for the Laplace-Beltrami operator
on $\S^{k-1}$, solving 
\begin{equation}
\label{eq:phi}
\begin{cases}
-\Delta_{\sigma}\phi=\mu\phi\\
\phi\in  H^1_{\mathcal G_k}(\S^{k-1})
\end{cases}~\!,
\quad \fint\limits_{\S^{k-1}}|{\phi}|^2~\!d\sigma=1~,
\quad\fint\limits_{\S^{k-1}}{\phi}~\!d\sigma=0,
\end{equation}
for some eigenvalue $\mu>0$, where $H^1_{{\mathcal G_k}}(\S^{k-1})$  is
the space of ${\mathcal G_k}$-invariant functions in $H^1(\S^{k-1})$.  
In particular one can take
\begin{equation}
\label{eq:first}
\mu=\inf~\left\{\dfrac{\int\limits_{\S^{k-1}}|\nabla_{\!\sigma} \f|^2d\sigma}{\int\limits_{\S^{k-1}}|\f|^2d\sigma}~~:~~
\f\in H^1_{{\mathcal G_k}}(\S^{k-1})~,~~\fint\limits_{\S^{k-1}}\f d\sigma=0\right\}~\!.
\end{equation}
To shorten notation we put 
$$
\phi_j(x)=\phi\big(\frac{x_j}{|x_j|}\big)~\quad \text{for $j=1,\dots, m$~,~~$x=(x_1,\dots, x_m)\ni \R^k_1\times\dots\R^k_m$.}
$$

\begin{Lemma}
\label{L:Lbs}
Let $\R^n=(\R^k)^m$, with $k\ge 2$, $m\ge 1$. Assume that $u\in\mathcal D^s$ is radially symmetric in each variable $x_j$, that is, $u(x_1,\dots,x_m)=u(|x_1|,\dots,|x_m|)$.  
The function
$$
\tu(x)= 
u(x)\sum_{j=1}^m \frac{|x_j|}{|x|}\phi_j(x) 
$$
belongs to  $\mathcal D^s$
and satisfies
\begin{gather}
\label{eq:normal}
\irn |x|^{-bq}|u|^{q-2}u\tu~\!dx=~0~,\quad\
\irn |x|^{-bq}|u|^{q-2}|\tu|^2~\!dx=\irn |x|^{-bq}|u|^{q}~\!dx~\!,
\\
\label{eq:uw}
\irn|\Dshalf \tu|^2~\!dx - \irn|\Dshalf u|^2~\!dx \le 
c_\mu\irn |x|^{-2s}|u|^2~\!dx~\!,
\end{gather}
where $q\in[2,\sstar]$, $b=\frac{n}{q}-\frac{n}{2}+s$ and the constant $c_\mu$  does not depend 
on $u$.
\end{Lemma}

\proof
We start by pointing out  the orthogonality relation
\begin{equation}
\label{eq:orthogonality2}
\irn \!\Big(\sum_{j,h=1}^{m}f_jg_h\phi_j\phi_h\Big) V~\!dx 
=\irn \!\Big(\sum_{j=1}^{m}f_jg_j\Big)V~\!dx,
\end{equation}
that holds for functions $V, f_j, g_
h$, each of them satisfying suitable summability assumptions and being radially symmetric in each variable $x_j$, $j=1,\dots,m$. To prove (\ref{eq:orthogonality2}) we first notice that  
$$
\sum_{h=1}^m ~\!f_j g_h\!\fint\limits_{\S_j^{k-1}}\phi_j\phi_h~\!d\sigma_j
=\sum_{h=1}^mf_j g_h~\!\delta_{jh}=f_jg_j\quad\text{for any $j=1,\dots, m$,}
$$
compare with (\ref{eq:phi}).
Since $Vf_jg_h$ is radially symmetric in $x_j$, we infer that
$$
\sum_{h=1}^m~\!\irn Vf_jg_h\phi_j\phi_h~\!dx=
\irn Vf_jg_j~\!dx~\!,
$$
so that (\ref{eq:orthogonality2}) follows by taking the sum for $j=1,\dots m$.

We are now in position to prove the lemma. The first equality in (\ref{eq:normal}) is immediate, because
for any 
index $j=1,\dots,m$ the function
$\phi_j=\phi(\sigma_j)$ has null mean on $\S^{k-1}_j$,
while $|x_j|^{-bq-1}|u|^q|x_j|$ is radially symmetric in
the variable $x_j\in\R^k_j$.

The second equality in (\ref{eq:normal}) follows from (\ref{eq:orthogonality2}). In fact,   $\sum_j|x_j|^2=|x|^2$ and thus
$$
\begin{aligned}
\irn |x|^{-bq}|u|^{q-2}|\tu|^2~\!dx=~\irn |x|^{-bq}|u|^{q}|x|^{-2}\big(\sum_{j,h=1}^m|x_j||x_h|\phi_j\phi_h\big)~\!dx=
\irn |x|^{-bq}|u|^q~\!dx.
\end{aligned}
$$

Next, let $w_u=w_u(x,y)$ be the Caffarelli-Silvestre extension of $u$. 
Since $w_u$ is uniquely determined
as the solution of a convex minimization problem, then  clearly $w_u(x,y)=w_u(|x_1|,\dots,|x_m|,y)$ 
for any $y>0$.

We introduce the following extension $\widetilde w$ of $\tu$,
$$
\widetilde w(x,y)=w_u(x,y)~\!\sum_{j=1}^m \frac{|x_j|}{\sqrt{|x|^2+y^2}}~\!\phi_j(x)~\!.
$$
From now on we simply write $w$ instead of $w_u$. 
It is also convenient to put
$$
\zeta=(x,y)\in \R^{n+1}_+~,\qquad  f_j(\zeta)=|\zeta|^{-1}|x_j|\ \ \text{for} \ \ j=1,\dots,m, \qquad F(\zeta)=\sum_{j=1}^{m}f_j\phi_j~\!.
$$
We claim that 
\begin{equation}
\label{eq:uw_better}
\iirn y^{1-2s}\big(|\nabla\widetilde w|^2-|\nabla w|^2\big)~\!dxdy\le (m\mu+m+1-2s)\iirn y^{1-2s}|\zeta|^{-2}|w|^2~\!dxdy~\!.
\end{equation}

To prove (\ref{eq:uw_better}) we use  (\ref{eq:orthogonality2}), (\ref{eq:dis}) and notice that  $\sum_{j}f_j^2=|\zeta|^{-2}|x|^2\le 1$
to get
\begin{gather}
\label{eq:F=1}
\irn y^{1-2s} F^2|\nabla w|^2~\!dx 
\le \irn y^{1-2s}|\nabla w|^2~\!dx;\\
\label{eq:partial}
-\irn F\partial_y F|w|^2~\!dx =
y\irn |w|^2~|\zeta|^{-4}|x|^2~\!dx\le \hat c\, y \irn |u|^2~|\zeta|^{-2}~\!dx.
\end{gather}
Since $\widetilde w=Fw$, from (\ref{eq:F=1}) we easily infer
$$
\iirn y^{1-2s}|\nabla\widetilde w|^2dxdy\le 
\iirn y^{1-2s}|\nabla w|^2dxdy
+\iirn y^{1-2s}\big(|w|^2 |\nabla\!F|^2+F\nabla F\cdot\nabla |w|^2\big)dxdy.
$$
Since (\ref{eq:partial}) and the Lebesgue dominated convergence theorem give 
\begin{gather*}
\lim_{y\to 0^+} \Big| \irn y^{1-2s}F\partial_y F |w|^2dx\Big|\le 
\hat c~\!\lim\limits_{y\to0^+} \irn y^{2-2s}\frac{|u|^2}{|x|^2+y^2}~\!dx=0\\
\lim_{y\to \infty} \Big| \irn y^{1-2s}F\partial_y F |w|^2dx\Big|\le 
\hat c~\!\lim\limits_{y\to\infty} \irn y^{2-2s}\frac{|u|^2}{|x|^2+y^2}~\!dx=0~\!,
\end{gather*}
(the summable majorant is $\frac{|u|^2}{|x|^{2s}}$),
we can integrate by parts on $\R^{n+1}_+$ to obtain 
$$
\begin{aligned}
\iirn y^{1-2s}~\!\big(|\nabla\widetilde w|^2-|\nabla w|^2\big)~\!dxdy
\le &~
\iirn |w|^2\Big(y^{1-2s}|\nabla\!F|^2-\div\big(y^{1-2s}F~\!\nabla\!F\big)\Big)~\!dxdy\\
=&~-\iirn |w|^2 F~\!\div\big(y^{1-2s}\nabla\!F\big)~\!dxdy~\!.
\end{aligned}
$$
To go further we compute 
$$
-\div\big(y^{1-2s}\nabla\!F\big)=y^{1-2s}
\sum_{j=1}^m\big(-\Delta(f_j\phi_j)-(1-2s)y^{-1}\partial_y(f_j\phi_j)\big)=
\sum_{j=1}^m g_j\phi_j~\!,
$$
where
$$
\begin{aligned}
g_j=&-\Delta f_j-(1-2s)y^{-1}\partial_y f_j+\mu|x_j|^{-2}f_j\\
=&~ |\zeta|^{-3} |x_j|^{-1}\big((n+1-2s)|x_j|^2+(\mu-k+1)|\zeta|^2\big)~\!.
\end{aligned}
$$
Now (\ref{eq:orthogonality2}) gives
\begin{equation*}
\label{eq:almost}
\iirn y^{1-2s}\big(|\nabla\widetilde w|^2-|\nabla w|^2\big)~\!dxdy\le
\iirn y^{1-2s}|w|^2\big(\sum_{j=1}^{m} f_jg_j\big)~\!dxdy.
\end{equation*}
Since
$$
\sum_{j=1}^{m} f_jg_j= \sum_{j=1}^{m} \big((n+1-2s)|x_j|^2+(\mu-k+1)|\zeta|^2\big)|\zeta|^{-4}
\le (m\mu+m+1-2s)|\zeta|^{-2},
$$
we readily obtain (\ref{eq:uw_better}).

By Lemma \ref{L:crucial}, inequality (\ref{eq:uw_better}) gives $\widetilde w\in\mathcal W^s$, thus $\tu=\widetilde w(~\!\cdot~\!,0)\in\mathcal D^s$. To conclude the proof, we compare 
the left-hand side of (\ref{eq:uw_better}) with (\ref{eq:CS_problem}) (with $\tu$ instead of $u$) and (\ref{quad_D}), and estimate the right-hand side by (\ref{eq:dis}).
\QED

In case $m=1$ (hence, $k=n$), we  have the following immediate corollary

\begin{Corollary}
\label{C:Lbs}
Assume $n\ge 2$ and  let $u\in \mathcal D^s$ be radially symmetric. Let ${\f_{\!\mu}}\in H^1(\S^{n-1})$ be a  nonconstant eigenfunction of the Laplace-Beltrami operator on $\S^{n-1}$ 
relative to the eigenvalue $\mu>0$ and  satisfying (\ref{eq:normalization}). Then the function $u{\f_{\!\mu}}\in \mathcal D^s$  satisfies (\ref{eq:uw_int}),
where $c_\mu$  does not depend 
on $u$.
\end{Corollary}

\section{Variational tools}
\label{S:variational}

We write $\mathcal G\prec O(n)$ if $\mathcal G$ is a closed subgroup of the orthogonal group  in $\R^n$
and put 
\begin{equation}
\label{eq:min_G}
\mathcal D^s_{\mathcal G}=\{u\in \mathcal D^s~|~ u\circ G\equiv u~~\text{for any} ~~G\in \mathcal G~\}~,\quad
S^{\mathcal G,\lambda}_q=\inf_{u\in\mathcal D^s_{\mathcal G}\atop u\neq 0} J_{\lambda}(u).
\end{equation}
One finds the
larger space $\mathcal D^s$ and the smallest constant $S^\lambda_q$ by choosing the trivial group
$\mathcal G$. 
The space $\mathcal D^s_{\text rad}$ of radial functions in $\mathcal D^s$ and the  infimum
$S^{{\text{rad}},\lambda}_q$ are recovered by taking $\mathcal G=O(n)$. 
Trivially, one has $\mathcal D^s_{\text rad}\subseteq \mathcal D^s_\mathcal G$ and $S^\lambda_q\le S^{\mathcal G,\lambda}_q\le S^{{\text{rad}},\lambda}_q$, for any $\mathcal G\prec O(n)$.

\begin{Remark}
\label{R:inv}
 Notice that $u\in \mathcal D^s_\mathcal G$ if and only if its Caffarelli-Silvestre extension $w_u(~\!\cdot~\!,y)$ is invariant with respect to the action of the group 
 $\mathcal G$, for any $y>0$.
\end{Remark}

\begin{Lemma}
\label{L:invariance}
Assume that (\ref{eq:bq}) is satisfied. Let $\lambda>-H_{\!s}$ and let $\mathcal G\prec O(n)$. 
Then the infimum $S^{\mathcal G,\lambda}_q$ is positive and attained. Moreover, if
$u\in \mathcal D^s_{\mathcal G}$ achieves $S^{\mathcal G,\lambda}_q$ then,
up to a Lagrange multiplier, $u$ is nonnegative and solves (\ref{eq:CKN}). 
Finally, for any $\tilde u\in \mathcal D^s_{\mathcal G}$ it holds that  
\begin{align}
\nonumber
(q-1)\frac{Q_\lambda(u)}{\||x|^{-b}u\|_q^q}\,\irn|x|^{-bq}|u|^{q-2}|\tilde u|^2~\!dx \,&\le
Q_\lambda(\tilde u)\\
&+(q-2)\frac{Q_\lambda(u)}{\||x|^{-b}u\|_q^{2q}}~\!\Big(\irn|x|^{-b}|u|^{q-2}u\tilde u~\!dx\Big)^2~\!,
\label{eq:J''}
\end{align}
where
$$
Q_\lambda(u):=\|\Dshalf u\|_2^2+\lambda\||x|^{-s}u\|_2^2~\!.
$$
\end{Lemma}

\proof
We already know that $S^{\mathcal G,\lambda}_q\ge S^\lambda_q>0$.
To show that the noncompact minimization problem in (\ref{eq:min_G}) 
has a solution we follow the outline of the proof of Theorem 0.1 in \cite{GM}, see also \cite{MNcn}.
By a standard convexity argument, we only need to construct a  minimizing sequence
that weakly converges to a nontrivial limit. 

We take a small number $\eps_0$ such that
\begin{equation}
\label{eq:eps0}
0<\eps_0<\frac12 {S^{\mathcal G,\lambda}_q}.
 \end{equation}
Hereafter, we denote by $B_R\subset \R^n$ the open ball of radius $R>0$ about the origin.
Since the ratio in (\ref{eq:min_problem}) is invariant with respect to the transforms 
$u(x)\mapsto  \alpha u(\beta x)$ (with $\alpha\neq 0, \beta>0$) of the space ${\mathcal D}^s_{\mathcal G}$ onto itself, we can find
a  minimizing sequence $u_h$ for $S^{\mathcal G,\lambda}_q$ such that
\begin{equation}
\label{eq:normaliz}
\||x|^{-b}u_h\|^{q}_q=({S^{\mathcal G,\lambda}_q)^\frac{q}{q-2}}~,~\qquad
Q_\lambda(u_h)=({S^{\mathcal G,\lambda}_q)^\frac{q}{q-2}}+o(1),  
 \end{equation}
\begin{equation}
\label{eq:Rda_sopra}
\eps_0^{\frac{q}{q-2}}\le \int\limits_{B_2}|x|^{-bq}|u_h|^q~\!dx
\le  \left(2\eps_0\right)^{\frac{q}{q-2}}\,, 
\end{equation}
and $u_h\to u$ weakly in  ${{\mathcal D}^s}$ for some $u\in {{\mathcal D}^s}$.
We only have to prove that $u\neq 0$.

We argue by contradiction. If $u_h\to 0$ weakly in ${\mathcal D^s_\mathcal G}$, we
can use Rellich theorem to get
that 
$|x|^{-b}u_h\to 0$ strongly in $L^q_{\rm loc}(\R^n\setminus\{0\})$. 
So, (\ref{eq:Rda_sopra}) implies
\begin{equation}
\label{eq:Rda_sopra2}
\eps_0^{\frac{q}{q-2}}\le \int\limits_{B_1}|x|^{-bq}|u_h|^q~\!dx +o(1).
\end{equation}
By Ekeland's variational principle we can assume that 
\begin{equation}
\label{eq:Rtu_equation}
\Ds u_h+\lambda|x|^{-2s}u_h-
|x|^{-bq}|u_h|^{{q}-2}u_h\to0\qquad\textrm{in ~${(\mathcal D^s_\mathcal G)}'$}.
\end{equation}
Take a radial function $\f\in {\cal C}^\infty_0(B_2)$ such that 
$\f\equiv 1$ on $B_1$. Then $\f^2u_h$ is a bounded sequence in $\mathcal D^s_\mathcal G$.
By  \cite[Lemma 2.1]{MNcn}, there exists a constant $c>0$ depending on $\f$ but not on $h$,
such that 
$$
\big|\langle \Ds \f u_h,\f u_h\rangle-\langle \Ds u_h,\f^2 u_h\rangle\big|^2\le  c~\!\langle \Ds u_h,u_h\rangle
\|u_h\|^2_{L^2(B_2)}.
$$
In particular, thanks to Rellich theorem we obtain
$$
\langle\Ds u_h,\f^2u_h\rangle =\langle\Ds (\f u_h),\f u_h\rangle+o(1)=\|\Dshalf(\f u_h)\|_2^2+o(1),
$$
that compared with the definition of $S^{\mathcal G,\lambda}_q$ leads to
\begin{equation}
\label{eq:numero}
\langle\Ds u_h,\f^2u_h\rangle+\lambda\| |x|^{-s}\f u_h\|_2^2
\ge {S^{\mathcal G,\lambda}_q}\||x|^{-b}\f u_h\|^2_q+o(1).
\end{equation}
On the other hand,  (\ref{eq:Rtu_equation}), H\"older's inequality and (\ref{eq:Rda_sopra}) give 
$$
\begin{aligned}
\langle\Ds u_h,\f^2u_h\rangle+\,&\,\lambda\| |x|^{-s}\f u_h\|_2^2
=\int\limits_{\R^n} |x|^{-bq}|u_h|^{{q}-2}|\f u_h|^2~\!dx+o(1)\\
\le&\, \Big(\int\limits_{~\!B_2}|x|^{-bq} |u_h|^q~\!dx\Big)^{\!\frac{q-2}{q}}
\||x|^{-b}\f u_h\|^2_q\le  2\eps_0 \||x|^{-b}\f u_h\|^2_q+o(1)~\!.
\end{aligned}
$$
Taking (\ref{eq:numero}) into account, we see that 
$$
{S^{\mathcal G,\lambda}_q}\||x|^{-b}\f u_h\|^2_q\le 2\eps_0 \||x|^{-b}\f u_h\|^2_q+o(1),
$$
which implies $\||x|^{-b}\f u_h\|_q=o(1)$ by (\ref{eq:eps0}).  
Since $\f\equiv 1$ on $B_1$, we infer from (\ref{eq:Rda_sopra2})
$$
\eps_0^{\frac{q}{q-2}}\le \int\limits_{B_1}|x|^{-bq}|u_h|^q~\!dx +o(1)\le \irn|x|^{-bq}|\f u_h|^q~\!dx +o(1)=o(1)~\!,
$$
a contradiction. Therefore, $S^{\mathcal G,\lambda}_q$ is achieved by some function $u\in \mathcal D^s_\mathcal G$.

Since $Q_\lambda(v)>Q_\lambda(|v|)$ for sign-changing function $v\in\mathcal D^s$, see \cite[Theorem 6]{MN-HSc}, we can assume that $u$  is {nonnegative}.

Next, for any $G\in O(n)$ we have $J_{\lambda}(u\circ G)=J_{\lambda}(u)$. As a consequence of the {\em Principle of symmetric criticality} \cite{P} we have that $\mathcal D^s_\mathcal G$ is 
a natural constraint for $J_{\lambda}$, and thus $u$ is a critical point for $J_{\lambda}$ on the whole $\mathcal D^s$. So, $u$ solves the fractional differential equation in (\ref{eq:CKN}), up to a Lagrange multiplier. 

Finally, the proof of (\ref{eq:J''}) is a simple computation, based on the fact that the function $f(t)=J_{\lambda}(u+t\tu)$ attains its minimum value at $t=0$, hence $f'(0)=0, f''(0)\ge 0$.
\QED

By taking $\mathcal G=\{\text{Id}_{\R^n}\}$ and then $\mathcal G=O(n)$ in Lemma \ref{L:invariance}, we immediately obtain the next existence result.

\begin{Corollary}
\label{C:existence}
Assume that (\ref{eq:bq}) is satisfied. If $\lambda>-H_{\!s}$, then the infimum
$S^\lambda_q$ in (\ref{eq:min_problem}) and the infimum
$$
S^{\text{rad},\lambda}_q=\inf_{u\in\mathcal D^s_{\text{rad}}(\R^n)\atop u\neq 0} J_{\lambda}(u)
$$
are attained by nonnegative solutions to (\ref{eq:CKN}).
\end{Corollary}

\begin{Remark}
Clearly Lemma \ref{L:invariance} and Corollary \ref{C:existence} hold also in case $n=1$, $0<s<\frac12$, 
with no changes in the proof.
\end{Remark}

In general,  solutions achieving the infima $S^{\mathcal G,\lambda}_q$ for different groups $\mathcal G\prec O(n)$ can coincide.
To obtain distinct solutions, we will use a special class of groups in $O(n)$. 

We need to introduce some notation. 
Recall that $n=km$. To any rotation $G\in O(k)$ and any permutation $P\in S_{m}$ we associate the rotation 
$$
PG\in O(n):\,(x_1,\dots,x_m)\mapsto P(Gx_1,\dots Gx_m), \qquad x_j\in\R^k_j.
$$
Further, to any subgroup $\mathcal G_k\prec O(k)$ we associate the following group of rotations in $\R^n$,
\begin{equation}
\label{eq:def}
\widetilde{\mathcal G}_{k}=\{PG~|~P\in  S_{m}~,~~ G\in \mathcal G_k~\}\prec O(n).
\end{equation}

Let $u\in \mathcal D^s_{\widetilde{\mathcal G}_k}$. Then $u$ is $\mathcal G_k$-invariant in each variable $x_j$
and is invariant with respect to any permutation of the $m$-tuple of vectors $(x_1,\dots, x_m)$. 

We will say that functions in $\mathcal D^s_{\widetilde{O}(k)}$ are 
{\em $k$-radially symmetric} (in \cite{N04} they are called $(k,0)$-radial, see also \cite{KN, NN}). 
Notice that a $k$-radially symmetric function depends only on $|x_1|$, \dots, $|x_m|$ and is invariant
with respect to permutations of its variables.

Clearly, $n$-radially symmetric functions are radial.

\begin{Lemma}[Basic Lemma]
\label{L:Cnk}
Let $n=km$ with $k\ge 2$ and $m\ge 1$, and let  ${\mathcal G_k}$ be a proper closed subgroup of $O(k)$. If $\lambda>0$ is large enough, then no extremal for 
$S^{\widetilde{\mathcal G}_{k},\lambda}_q$ can be $k$-radially symmetric.
\end{Lemma}

\proof
Let $u\in\mathcal D^s_{\widetilde{\mathcal G}_k}$ be a $k$-radially symmetric function
achieving the best constant $S^{\widetilde{\mathcal G}_{k},\lambda}_q$ for some $\lambda>0$. 
We take $\tu$ as in Lemma \ref{L:Lbs}, where $\phi$ is given by (\ref{eq:phi}) while $\mu$ is defined in (\ref{eq:first}). Clearly $\tu\in\mathcal D^s_{\widetilde{\mathcal G}_k}$, 
so that formula (\ref{eq:J''}) in Lemma \ref{L:invariance} applies.

Using (\ref{eq:normal}) we can rewrite (\ref{eq:J''}) as follows:
$$
(q-1)Q_\lambda(u)\le Q_\lambda(\tu)
\le Q_\lambda(u)+c_\mu \irn |x|^{-2s}|u|^2~\!dx
$$
(the last inequality follows by (\ref{eq:uw})), i.e. 
$$
(q-2)Q_\lambda(u)\le c_\mu \irn |x|^{-2s}|u|^2~\!dx~\!.
$$
Now we use the Hardy inequality to estimate
$Q_\lambda(u)> (H_{\!s}+\lambda)\||x|^{-s}u\|_2^2$. We infer that 
$(q-2)(H_{\!s}+\lambda)<c_\mu$. We proved that if $\lambda\ge- H_{\!s}+\frac{c_\mu}{q-2}$ then
a $k$-radially symmetric function cannot provide the constant $S^{\widetilde{\mathcal G}_{k},\lambda}_q$.
\QED

\section{Proofs of the main results}
\label{S:proofs}

{\bf Proof of Theorem \ref{T:BS}.}
This is an immediate consequence of Lemma \ref{L:Cnk}, with
$m=1$, $k=n$ (so that $k$-radially symmetric functions are radially symmetric) and $\mathcal G=\{\text{Id}_{\R^n}\}$. 
\QED

\medskip
\noindent
{\bf Proof of Theorem \ref{T:hradial}.}
We identify $\R^n$ with $(\R^2)^m$, $m\ge 1$ and introduce the polar coordinates $x_j=(r_j,\sigma_j)$ for points in $\R^2$.

Following \cite{N04}, for any integer $t> 1$ we denote by 
$\mathcal G_2^t\prec O(2)$ the group generated by a rotation of $\frac{2\pi}t$ and by $\widetilde{\mathcal G}_2^t$ the corresponding subgroup of $O(n)$, see (\ref{eq:def}),
Then we denote by $u_t$ the nonnegative solution to (\ref{eq:CKN})
 solving the minimization problem
$$
J_\lambda(u_t)= \min_{u\in\mathcal D^s_{\widetilde{\mathcal G}_2^t}
} J_{\lambda}(u)~\!,
$$
compare with Lemma \ref{L:invariance}. We prove that for any pair of distinct integers $t, T$, 
the functions $u_T$, $u_t$ cannot coincide up to rotations, provided that $\lambda$ is large enough.

First, we face the case when $T=ht$, $h>1$, is a multiple of the integer $t$. Recall Lemma \ref{L:Cnk} and find $\lambda>0$ large, 
so that $u_T$ is not $k$-radially symmetric. We introduce the function
$$
v_t\big((r_1,\sigma_1),\dots , (r_m,\sigma_m)\big)=u_T\big((r_1,\frac {1}{h}\sigma_1),\dots , (r_m,\frac{1}{h}\sigma_m)\big).
$$
Easily, $v_t\in \mathcal D^s_{\widetilde{\mathcal G}_2^t}$ and
$$
\irn |x|^{-2s}|v_t|^2~\!dx=\irn |x|^{-2s}|u_T|^2~\!dx~,\quad
\irn |x|^{-bq}|v_t|^q~\!dx=\irn |x|^{-bq}|u_T|^q~\!dx~\!.
$$ 

Let $w_T=w_{u_T}$ be the Caffarelli-Silvestre extension of $u_T$ and consider the function
$$
\widetilde w_{t}\big((r_1,\sigma_1),\dots , (r_m,\sigma_m),y\big)=
w_T\big((r_1,\frac {1}{h}\sigma_1),\dots , (r_m,\frac{1}{h}\sigma_m),y\big)~\!. 
$$
We have $\widetilde w_{t}(x,0)=v_t(x)$. Moreover, from the formula
\begin{equation*}
|\nabla w|^2 = \sum\limits_{j=1}^m \Big(|\partial_{r_j}w|^2 + \frac{1}{r_j^2} |\partial_{\sigma_j}w|^2\Big)+|\partial_yw|^2
\end{equation*}
and taking into account the relation $\partial_{\sigma_j}\widetilde w_{t}=\frac 1h\partial_{\sigma_j}w_T$, we get 
\begin{equation}
\label{eq:wt}
\iirn y^{1-2s}|\nabla \widetilde w_{t}|^2\,dxdy<\iirn y^{1-2s}|\nabla  w_T|^2\,dxdy,
\end{equation}
because $w_T(~\!\cdot~\!,y)$ cannot be $k$-radially symmetric, and thus $\partial_{\sigma_j}w_T\not\equiv0$.

Formulae (\ref{eq:CS_problem}) with $\tu_t$ instead of $u$, (\ref{eq:wt}) and (\ref{quad_D}) give
\begin{align*}
\irn|\Dshalf v_t|^2~\!dx \le &\ C_s\iirn y^{1-2s}|\nabla \widetilde w_{t}|^2~\!dxdy \\
< &\ C_s\iirn y^{1-2s}|\nabla  w_T|^2~\!dxdy =\irn|\Dshalf u_T|^2~\!dx~\!,
\end{align*}
and we infer that 
\begin{equation}
\label{eq:sunday}
J_\lambda(u_t)\le J_\lambda(v_t)< J_\lambda(u_T)
\end{equation}
for $\lambda$ large enough. Thus, the statement is proved for $T=ht$.

In the case of general distinct intergers $T,t$ we define $\hat t$ as the least common multiple of the pair $T, t$. If $u_{T}=u_t\circ G$ for some rotation $G\in O(n)$ 
then $u_t$ is $\widetilde{\mathcal G}_2^{\hat{t}}$-invariant. But this implies $J_\lambda(u_{\hat{t}})\le J_\lambda(u_t)=J_\lambda(u_{T})$, that is impossible for $\lambda$ 
large enough by (\ref{eq:sunday}). The proof is complete.
\QED

Before the proof of Theorem \ref{T:dimension}, we make a remark. 

\begin{Remark}
Notice that if $n$ has $\ell$ distinct divisors then
it is easy to obtain  $\ell$ distinct solutions for $\lambda$ large enough. 
Let $1<k_j<n$, $j=1,\dots, \ell-2$, be the distinct nontrivial divisors of $n$, so that $\R^n=(\R^{k_j})^\frac{n}{k_j}$ . 
Thanks to Lemmata \ref{L:invariance} and \ref{L:Cnk}, for $\lambda$ large 
the $\ell$ solutions achieving the best constants $S^{{\mathcal G},\lambda}_q$ for $\mathcal G=\{\text{Id}\}$, $\mathcal G=\widetilde {O}(k_j)$ and $\mathcal G=O(n)$
are distinct modulo rotations.
\end{Remark}

To manage a general case, including prime dimensions, some extra argument is needed.

\medskip
\noindent
{\bf Proof of Theorem \ref{T:dimension}.}
Take any $u\in \mathcal D^s$ and its Caffarelli-Silvestre extension $w_u$. 
We denote by $u^*$ and $w^*(~\!\cdot~\!,y)$ the
symmetrization along spheres of $u$ and $w_u(~\!\cdot~\!,y)$, respectively.
By \cite[Theorem 2.31]{Kawohl}, such symmetrization diminishes the $L^2(\R^n)$-norm of $\nabla_{\!x} w$. Therefore
we get
$$
\iirn y^{1-2s}|\nabla w^*|^2~\!dxdy\le\iirn y^{1-2s}|\nabla w_u|^2~\!dxdy.
$$
Since evidently the spherical symmetrization keeps weighted norms,
$$
\irn |x|^{-2s}|u^*|^2~\!dx=\irn |x|^{-2s}|u|^2~\!dx~,\quad
\irn |x|^{-bq}|u^*|^q~\!dx=\irn |x|^{-bq}|u|^q~\!dx~\!,
$$ 
formulae (\ref{eq:CS_problem}) and (\ref{quad_D}) imply $J_\lambda(u^*)\le J_\lambda(u)$.
Therefore, if $U$ achieves $S^{\lambda}_q$ then it is axisymmetric, i.e., up to rotations, $U(x)=U(|x|,\vartheta)$ where $\vartheta=\cos^{-1}(x_n/|x|)$ is the angle between $x$ and the axis $Ox_n$. 
Moreover, if $U$ is not radial (that holds for $\lambda$ large enough) then it is strictly monotone with respect to $\vartheta$.

Now we define four groups such that corresponding minimizers are different for $\lambda$ large enough in arbitrary dimension. Besides full group $O(n)$ and the trivial group, they are 
$O(n-1)\times \mathbb{Z}_2$ (modulo rotations, corresponding functions depending only on $|x|$ and $|\cos(\vartheta)|$) and the symmetry group of the right simplex. 

By Lemma \ref{L:Cnk}, minimizers corresponding to the last three groups cannot be radial for $\lambda$ large enough. Since the global minimizer is axisymmetric and monotone with respect to $\vartheta$, 
it cannot be invariant with respect to the last two groups. Similarly, the minimizer generated by $O(n-1)\times \mathbb {Z}_2$ is monotone with respect to $\vartheta$ in both half-spaces, and so it 
cannot be invariant with respect to the simplex group.
\QED

\subsection*{Conclusions}
We furnished a powerful tool that can be applied to a large class of variational equations driven by the fractional Laplace operator $\Ds$
of order $s\in(0,1)$. One of the main steps is Lemma \ref{L:Lbs}, that allows to compare the $L^2$ norms of
$\Dshalf u$ and $\Dshalf (u\f_\mu)$, where $u$ is a given function of $k$ variables $x_j\in\R^m$ depending only on $|x_j|$,
and $\f_\mu$ is a suitably normalized eigenfunction of the Laplace-Beltrami operator on the unit sphere in $\R^{mk}$.

Then we take as model a nonlocal variational equation related to the fractional Caffarelli-Kohn-Nirenberg inequality 
to illustrate how Lemma \ref{L:Lbs} can be used in order to obtain symmetry breaking and multiplicity phenomena.

Differently from the local case $s=1$, an efficient Emden-Fowler transform is not available and positive radial solutions to (\ref{eq:CKN}) are not explicitly known;
their uniqueness (up to dilations) is an open question as well, that makes the problem more challenging.

\subsection*{Acknowledgements}

R.M. was partially supported by PRID project {\em VARPROGE}.

This work does not have any conflicts of interest

\footnotesize
\label{References}

\end{document}